\documentclass[12pt,reqno]{article} 
\usepackage{amssymb,amsmath,amsthm,graphicx}
\usepackage{color}
\usepackage{geometry}
 \newcommand{\bfR}{{\bf R}}
 \newcommand{\bfN}{{\bf N}}
 \newcommand{\bfZ}{{\bf Z}}

 \newcommand{\la}{\lambda}

 \renewcommand{\th}{\theta}

 \newcommand{\Th}{\Theta}
 \newcommand{\sig}{\sigma}
 \newcommand{\bk}{{\bf k}}
\newcommand{\bp}{{\bf p}}
\newcommand{\bq}{{\bf q}}

\newcommand{\bv}{{\bf v}}
\newcommand{\by}{{\bf y}}

 \newcommand{\harth}{{\overset{\rightharpoonup}{\theta}}}
 
 \renewcommand{\d}{\displaystyle}
 \newcommand{\rme}{{\rm e}}

 \date{} 
\numberwithin{equation}{section}

\allowdisplaybreaks

\begin{document}

 \title{\bf Periodic spectrum of $n$-cubic quantum graphs}
 \author{Chun-Kong\ Law$^1$, Yu-Chun\ Luo$^2$ and Tui-En Wang$^3$}
 \maketitle
 \begin{center}
 \bf{Abstract}
 \end{center}
 \par We study the spectrum of some periodic differential operators, in particular the periodic Schr\"{o}dinger operator acting on infinite
 $n$-cubic graphs. Using Floquet-Bloch theory, we derive and analyze on the dispersion relations of the periodic quantum graph generated
 by 2-dimensional rectangles, and also $n$-cubes. Our proof is analytic. These dispersion relations define the spectra of the associated periodic operator, thus
 facilitating further analysis of the spectra.
 \\[0.5in]

 \footnote{$^1$Department of Applied Mathematics, National Sun
 Yat-sen University, Kaohsiung, Taiwan 80424, R.O.C. Email:
 law@math.nsysu.edu.tw}
 \footnote{$^2$Department of Applied Mathematics, National Sun
 Yat-sen University, Kaohsiung, Taiwan 80424, R.O.C. Email:
 leoredro@gmail.com}
 \footnote{$^3$Department of Applied Mathematics, National Sun
 Yat-sen University, Kaohsiung, Taiwan 80424, R.O.C. Email:
 miteric1024@gmail.com}
 \newpage

 \setcounter{equation}{0}
 \section{Introduction}
 \hskip0.25in
 Recently there have been a lot of studies on quantum graphs (see \cite{BK} and its reference). In particular, Kuchment and Post \cite{KP}
 employed this theory to understand the spectrum of Schr\"{o}dinger operators acting on graphene,
 where the material is represented by a periodic regular hexagonal graphs. Applying the Floquet-Bloch theory, they
  derived the dispersion relation for the periodic hexagonal quantum graphs, which is exactly the quantum network model for graphene.
  Their derivation goes through Hill's discriminant. Through a previous work, we have employed the idea that
  the eigenvalues of a compact quantum graph are zeros of the characteristic function \cite{LP}.  We found that
  this characteristic function approach is an efficient method to evaluate the spectrum of periodic quantum graphs.

  An Archimedean tiling (or uniform tiling) is a pattern of regular polygons with the same configuration that covers
  the whole plane and fits around every vertex of the polygons.  There are totally 11 of them \cite{GS}, including the square tiling and
  hexagonal tiling.  In \cite{LJL1,LJL2},  some of us used the characteristic function method to derive the dispersion relations for
  all the other nine Archimedean tilings, for example, the triangular tiling ($3^6$), elongated triangular tiling ($3^3,4^2$),
  trihexagonal tiling ($3,6,3,6$) and truncated square tiling ($4,8^2$),
  However, some of the matrices involved are too big (dim$=10,\ 12$, up to $36$).  We need to use the symbolic manipulation software Mathematica to evaluate
  their determinants and algebraic simplification.

 On the interval $[0,a]$, we let $C(x,\rho)$ and $S(x,\rho)$  {($\la=\rho^2$)} be the solutions of
 $$
 -y''+q y=\lambda y
 $$
 such that $C(0,\rho)=S'(0,\rho)=1$, $C'(0,\rho)=S(0,\rho)=0$. We call them cosine-like function and sine-like function
 respectively. They can also be expressed in an integral form:
 \begin{eqnarray}
 C(x,\rho)&=&\displaystyle\cos(\rho x)+\frac1{\rho}\int^x_0\sin(\rho (x-t))q(t)C(t, \rho)dt,\label{eq1.11}\\
 S(x, \rho)&=&\displaystyle \frac{\sin(\rho x)}{\rho}+\frac1{\rho}\int^x_0\sin(\rho (x-t))q(t)S(t, \rho)dt.
 \label{eq1.12}
 \end{eqnarray}
 Note if $q$ is even, then $ \d S'(a,\rho)=C(a,\rho)$
 by \cite[p.8]{MW} (see also \cite{PR13}).

  \newtheorem{th1.0}{Theorem}[section]
 \begin{th1.0}[\cite{LJL1}]
 \label{th1.0}
  Assume that all the $q_j$'s are identical to $q$, and even. We also let $\theta_1,\ \theta_2\in[-\pi,\pi]$.
 \begin{enumerate}
 \item[(a)] For the triangular tiling, the dispersion relation of the associated periodic quantum graph is
 $$
 S^2 \left(3S'+1-4 \cos(\frac{\th_1}{2})\cos(\frac{\th_2}{2})\cos(\frac{\th_2-\th_1}{2})\right)=0.
 $$
 \item[(b)] For the elongated triangular tiling, the dispersion relation of the associated periodic quantum graph is given by
  $$
  S^3\, \{25(S')^2-20\cos\theta_1S'-8\cos(\frac{\theta_1}{2})\cos(\frac{\theta_2}{2})\cos(\frac{\theta_1-\theta_2}{2})-4\cos^2\theta_1-1\}=0.
 $$
 \item[(c)] For the truncated square tiling, the dispersion relation of the associated periodic quantum graph is given by
 $$
 S^2\left\{ 81 S'^4-54 S'^2-12S'(\cos\th_1+\cos\th_2)+1-4\cos\th_1\, \cos\th_2\right\} = 0.
 $$
 \item[(d)] For trihexagonal tiling, the dispersion relation of the associated periodic quantum graph is given by
 $$
 S^3(2S'+1)\, \left(2 S'^2-S'-\cos\frac{\theta_2}{2}\cos\frac{\theta_1-\theta_2}{2})\right)=0.
 $$
 \end{enumerate}
 \end{th1.0}

  Since these dispersion relations define the spectra of the corresponding periodic differential operators, we can understand more about the spectra by analyzing
  on the variety of these polynomials of $S(a,\rho)$ and $S'(a,\rho)$. Hence the spectra of these periodic quantum graphs can be expressed in terms of the periodic spectrum
  on the line.    We remark that in \cite{LJL1}, more general dispersion relations were given for the above four Archimedean tiling, where
   the potential functions are not required to be identical or even. For example, the dispersion relation for the triangular tiling
 is
 \begin{eqnarray*}
 \lefteqn{ S_1'S_2S_3+S_1 S_2'S_3+S_1 S_2 S_3'+C_1 S_2S_3+C_2S_1S_3+C_3S_1S_2}\\
  &&-2S_1 S_2\, \cos(\th_1-\th_2)-2 S_2 S_3\, \cos\th_1-2 S_1 S_3\cos\th_2\  = \ 0.
  \end{eqnarray*}

  In this paper, we plan to employ the same method to study the periodic Schr\"odinger operator acting on $n$-cubic quantum graphs. We shall
  first study the rectangular graph ($n=2$), then the general $n$-cubic graph.  We do not need to require $q_i$'s to be identical and even.
  The edgelengths $a_i$ may be different too.
  (But when they are so, the dispersion relations become a lot simpler, as we shall see.) Furthermore our analysis is rigorous and general, \underline{with
  no resort} to any symbolic manipulation software. For we find that the determinant can be decomposed
  into two smaller determinants: one with nontrivial first column and diagonal, another with nontrivial last row and sub-diagonal. By Lemma~\ref{lem3.5}, we computed these
  special determinants and achieve our goal.

  Given an infinite graph $\Gamma=E(\Gamma)\cup V(\Gamma)$ generated by an Archimedean tiling, we let $H$ denote a Schr\"odinger operator on , i.e.,
 $$
 H \, \by(x)=-\frac{d^2}{d x^2}\by(x)+\bq(x)\, \by(x),
 $$
  {where $ \bq\in L^2_{loc}(\Gamma)$ is periodic on the tiling (explained below), and} the domain $D(H)$ consists of all  {admissible} functions $\by(x)$ (union of functions $y_e$ for each edge $e\in E(\Gamma)$) on $\Gamma$ in the sense that
 \begin{enumerate}
 \item[(i)] $\displaystyle   y_e\in {\cal H}^2(e)$ for all $e\in E(\Gamma)$;
 \item[(ii)] $\displaystyle   \sum_{ {e\in E(\Gamma)}} \| y_e\|^2_{{\cal H}^2(e)} <\infty$;
 \item[(iii)] Neumann vertex conditions (or continuity-Kirchhoff conditions at vertices), i.e., for any vertex $\bv\in V(\Gamma)$,
 $$
 y_{e_1}(\bv)=y_{e_2}(\bv)\qquad \mbox{ and }\qquad \sum_{ {e\in E(\Gamma)}} y_e'(\bv)=0.
 $$
 Here $y_e'$ denotes the directional derivative along the edge $e$ from $\bv$,  {and $\| \cdot\|^2_{{\cal H}^2(e)}$ denotes the Sobolev norm of 2 distribution derivatives.}
 \end{enumerate}
 The potential function $\bq$ is said to be periodic if
 $\bq(x+\bp\cdot \vec{\bk})=\bq(x)$ for all $\bp\in \bfZ^2$ and all $x\in \Gamma$.
  Take the quasi-momentum $\harth=(\theta_1,\theta_2)$ in the Brillouin zone
 $B=[-\pi,\pi]^2$. Let $H^{\harth}$ be the Bloch  {Hamiltonian} that acts on  {$L^2(\Gamma)$, and the dense domain $D(H^{\harth})$ consists of admissible functions $\by$ which satisfy the Floquet-Bloch condition}
 \begin{equation}
 \by(x+\bp\cdot \vec{\bk})= \rme^{i(\bp\cdot\harth)}\by(x),
 \label{eq1.02}
 \end{equation}
 for all $\bp\in \bfZ^2$ and all $x\in \Gamma$.  Such functions are uniquely determined by their restrictions on the fundamental domain $W$. Hence for fixed $\harth$,
 the operator $H^\harth$ has purely discrete spectrum $\sig(H^\harth)=\{\la_j(\harth):\ j\in \bfN\}$, where
 $$
  {\la_1(\harth)\leq \la_2(\harth)\leq \cdots \leq \la_j(\harth)\leq\cdots},\quad \mbox{ and } \la_j(\harth)\rightarrow\infty\mbox{ as } j\to\infty.
 $$
 
 By an analogous argument as in \cite[291]{RS}, $H$ is unitarily equivalent to a direct integral of $\displaystyle   H^\harth$, denoted as
 $$
 UHU^*=\int_B^\oplus H^{\harth}\, d\harth.
 $$
 for some unitary operator $U$.  Thus by \cite{BK,RS},
 $$
 \bigcup \{\sig(H^{\harth}):\ \harth\in [-\pi,\pi]^2\}=\sigma(UHU^*)=\sig(H).
 $$
 Furthermore,  it is known that singular continuous spectrum is absent in $\sig(H)$. So  $\sig(H)$ consists of only point spectrum and absolutely continuous spectrum \cite[Theorem 4.5.9]{K93}.

 Our main theorems are Theorem \ref{th3.2} and Theorem \ref{th3.4}, where the dispersion relations for rectangular graphs and $n$-cubic graphs are given.
 These materials occupy section 2 and section 3 respectively. We also study the point spectrum and absolutely continuous spectrum, with the help
 of these dispersion relations. Finally there is a section for concluding remarks.

 It has been known that molecules like graphene
 or boron nitride (BN) have crystal lattices \cite{ALM, LAM}.  By a
 quantum network model (QNM), the wave functions at the bonds (Bloch waves) satisfy a Schr\"{o}dinger equation along the lines, and
 continuity and Kirchhoff conditions at the vertices. So under this QNM, the spectrum represent the energy of the wave functions in
 the molecule. Thus the associated spectral analysis has physical significance in quantum mechanics. In \cite{ALM}, graphene is associated with hexagonal tiling, and
 the potential function is approximated by
 $$
 q(x)=-0.85+\displaystyle\frac{d}{1.34}\sin^2\left(\frac{\pi x}{d}\right)
 $$
 where $d$ is the distance between neighboring carbon atoms, and $d=1.43$ {\AA}.  Hence the spectra of their periodic graphs, as in the case of graphene, have physical meaning and their analyses worthwhile.  Moreover many crystal lattices are 3-dimensional, and the cube is the simplest tiling there.  So we believe our results have both mathematical and physical
 significance.

 We also remark that the Dirac point and other properties of some crystal lattices were also explored in \cite{K2013} and \cite{FW12}. We shall pursue on this issue later.

 \setcounter{equation}{0}
 \section{Dispersion relation for periodic rectangular graphs}
 \begin{figure}[h!]\label{fig3.1}
 \centering\includegraphics[width=6.2cm,height=6cm]{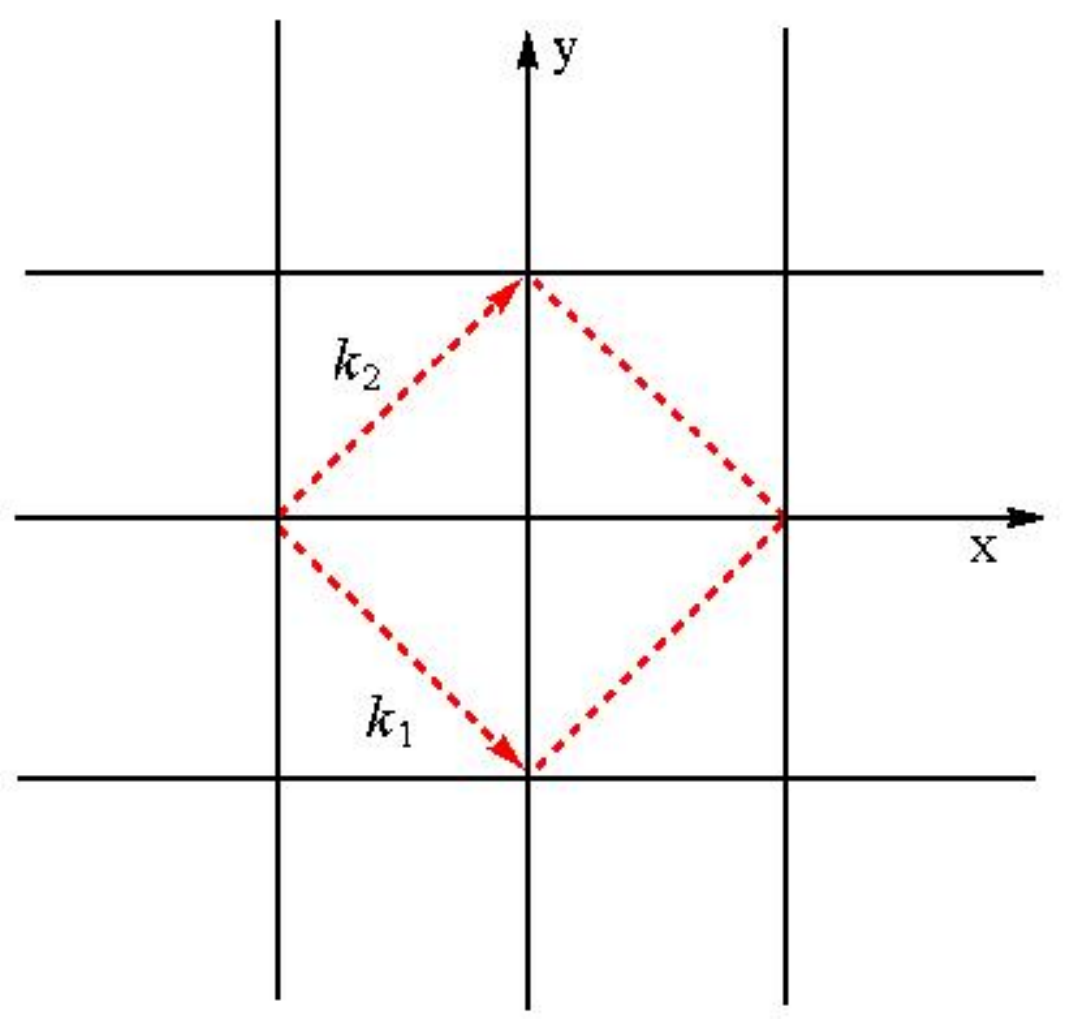}
 \caption{Fundamental domain of periodic rectangular graphs}
 \end{figure}

 \hskip0.25in Let $\Gamma$ be the infinite rectangular graph with edges $a_1$ and $a_2$ and $W$ be
 the fundamental domain of $\Gamma$ formed by two vectors $k_1=(a_1,-a_2)$ and $k_2=(a_1,a_2)$ such that
 $\Gamma=\displaystyle\bigcup_{\overset{\rightharpoonup}{p}\in\mathbb{Z}^2}\overset{\rightharpoonup}{p}\circ W$. Consider $Hy''(x):=-y''(x)+q(x)y(x)$.
 Let $\displaystyle H=\int^\oplus_{[-\pi, \pi]^2}H^{\overset{\rightharpoonup}{\theta}}, d\,\theta$ where in local coordinates,
 \begin{equation}
 H^{\overset{\rightharpoonup}{\theta}} y(x):=-y_j''(x)+q_j(x)y_j(x)=\lambda y_j(x,) \label{eq3.1}
 \end{equation}
 where $q_j(x)\in L^2(0, a_j),\ j=1,\ 2,\ 3,\ 4$, with the Floquet-Bloch conditions imposed:
 $$y_j(x+p_1k_1+p_2k_2)=e^{i(p_1\theta_1+p_2\theta_2)}y_j(x),$$
 where $\harth=(\theta_1, \theta_2)\in[-\pi,\pi]^2,\ (p_1,p_2)\in\mathbb{Z}^2$.

 \begin{figure}[h!]\label{fig3.2}
 \centering\includegraphics[width=9cm,height=6cm]{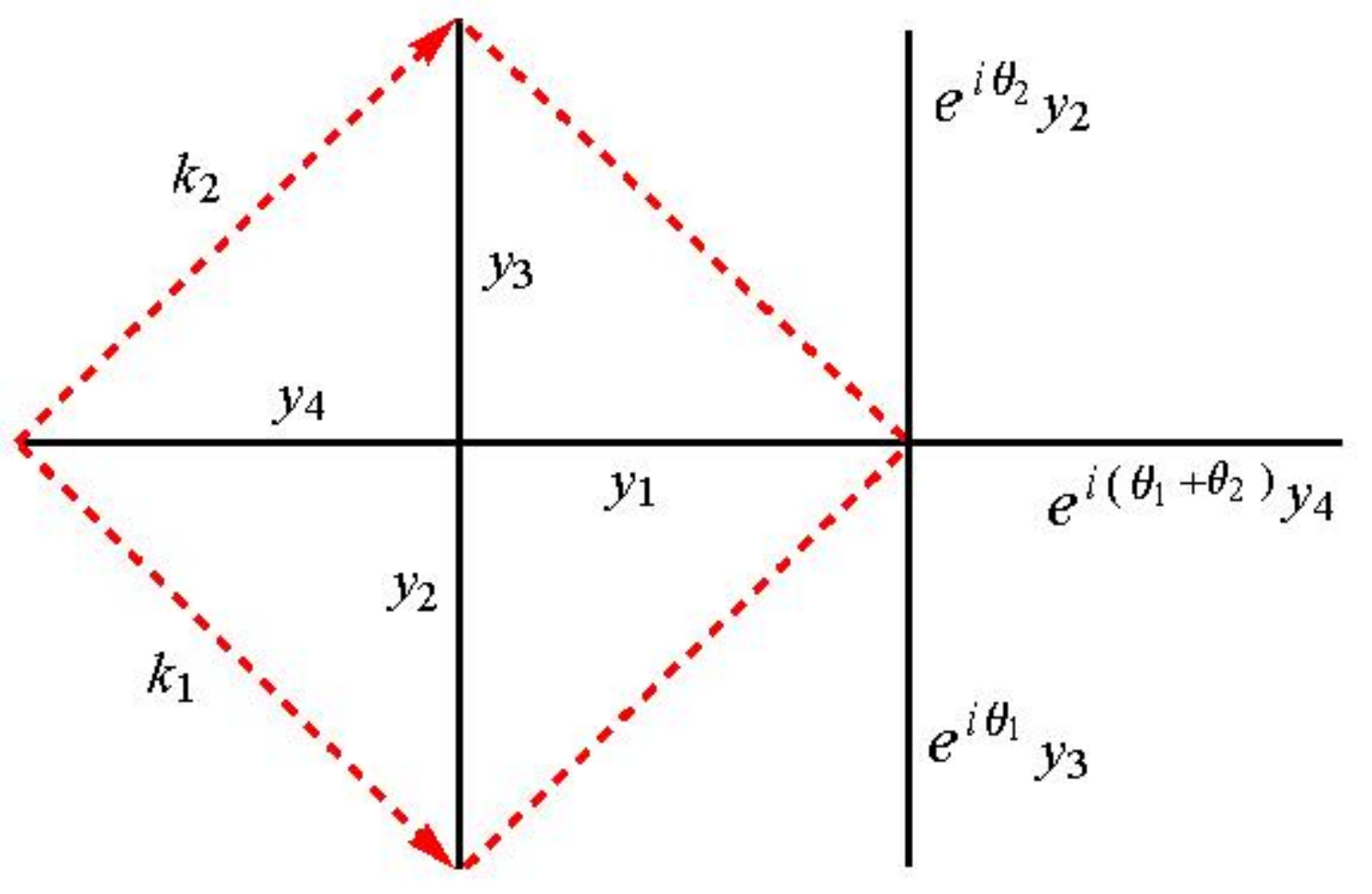}
 \caption{Floquet-Bloch conditions on periodic rectangular graphs}
 \end{figure}

 Coupled with the Neumann vertex conditions, we can make some analysis of the eigenfunction
 $y(\cdot,\lambda)$ in local coordinates. At ${\bf 0}$, we have
 \begin{equation}
 \left\{
 \begin{array}{l}
 y_1(0,\lambda)=y_2(0,\lambda)=y_3(0,\lambda)=y_4(0,\lambda)\\
 \displaystyle\sum^4_{j=1}y_j'(0,\lambda)=0
 \end{array}
 \right.
 \label{eq3.2}
 \end{equation}

 Then at $(a_1,0)$, by the Floquet-Bloch conditions, continuity and Kirchhoff conditions, we have
 \begin{equation}
 \left\{
 \begin{array}{l}
 y_1(a_1,\lambda)=e^{i\theta_2}y_2(a_2,\lambda)=e^{i\theta_1}y_3(a_2,\lambda)=e^{i(\theta_1+\theta_2)}y_4(a_1,\lambda)\\
 y_1'(a_1,\lambda)+e^{i\theta_2}y_2'(a_2,\lambda)+e^{i\theta_1}y_3'(a_2,\lambda)+e^{i(\theta_1+\theta_2)}y_4'(a_1,\lambda)=0
 \end{array}
 \right.
 \label{eq3.3}
 \end{equation}

 Let $\rho=\sqrt{\lambda}$, then any solution $y_j$ of (\ref{eq3.1}) has to satisfy:
 \begin{equation*}
 y_j(x, \rho)=A_j(\rho)C_j(x, \rho)+B_j(\rho)S_j(x, \rho),\ j=1,\ldots,\ 4
 \end{equation*}
 where the functions $S_j(x,\rho)$ and $C_j(x,\rho)$ satisfy \eqref{eq1.11}-\eqref{eq1.12}
 Since $C_j(0,\rho)=S_j'(0,\rho)=1$, $C_j'(0,\rho)=S_j(0,\rho)=0$, then we may rewrite (\ref{eq3.2}) and (\ref{eq3.3}) as follow:
 (Here, we write $C_j=C_j(a_j,\rho)$ and $S_j=S_j(a_j,\rho)$ for short.)
 \begin{equation}\label{eq3.4}
 \left\{
 \begin{array}{l}
 A_1(\rho)=A_2(\rho)=A_3(\rho)=A_4(\rho)\\
 \displaystyle\sum^{4}_{j=1}B_j(\rho)=0\\
 A_1(\rho)C_1+B_1(\rho)S_1=e^{i \theta_2}(A_2(\rho)C_2+B_2(\rho)S_2)=e^{i \theta_1}(A_2(\rho)C_3+B_2(\rho)S_3)\\
                          \qquad=e^{i (\theta_1+\theta_2)}(A_4(\rho)C_4+B_4(\rho)S_4)\\
 A_1(\rho)C_1'+B_1(\rho)S_1'+e^{i \theta_2}(A_2(\rho)C_2'+B_2(\rho)S_2')+e^{i \theta_1}(A_2(\rho)C_3'+B_2(\rho)S_3')\\
                          \qquad+e^{i (\theta_1+\theta_2)}(A_4(\rho)C_4'+B_4(\rho)S_4')=0
 \end{array}
 \right.
 \end{equation}

 Then we have the characteristic function $\Phi(\rho)$ as below:

 \begin{eqnarray*}
 \Phi(\rho)&=&
 \left|
 \begin{array}{cccccccc}
 1 & -1 & 0 & 0 & 0 & 0 & 0 & 0 \\
 1 & 0 & -1 & 0 & 0 & 0 & 0 & 0 \\
 1 & 0 & 0 & -1 & 0 & 0 & 0 & 0 \\
 0 & 0 & 0 & 0 & 1 & 1 & 1 & 1 \\
 C_1 & -e^{i \theta _2}C_2 & 0 & 0 & S_1 & -e^{i \theta _2}S_2 & 0 & 0 \\
 C_1 & 0 & -e^{i \theta _1}C_3 & 0 & S_1 & 0 & -e^{i \theta _1}S_3 & 0 \\
 C_1 & 0 & 0 & -e^{i (\theta_1+\theta_2)}C_4 & S_1 & 0 & 0 & -e^{i (\theta_1+\theta_2)}S_4 \\
 C_1' & e^{i \theta _2}C_2' & e^{i \theta _1}C_3' & e^{i (\theta_1+\theta_2)}C_4' & S_1' & e^{i \theta _2}S_2' & e^{i \theta _1}S_3' & e^{i
 (\theta_1+\theta_2)}S_4'
 \end{array}
 \right|\\
 &=&\left|
 \begin{array}{cccccccc}
 0 & 1 & 1 & 1 & 1 \\
 C_1-e^{i \theta _2}C_2 & S_1 & -e^{i \theta _2}S_2 & 0 & 0 \\
 C_1-e^{i \theta _1}C_3 & S_1 & 0 & -e^{i \theta _1}S_3 & 0 \\
 C_1-e^{i (\theta_1+\theta_2)}C_4 & S_1 & 0 & 0 & -e^{i (\theta_1+\theta_2)}S_4 \\
 C_1'+e^{i \theta _2}C_2'+e^{i \theta _1}C_3'+e^{i (\theta_1+\theta_2)}C_4' & S_1' & e^{i \theta _2}S_2' & e^{i \theta _1}S_3' & e^{i (\theta_1+\theta_2)}S_4'
 \end{array}
 \right|
 \end{eqnarray*}

 We use cofactor expansion for the first column and then along the last row. Observing that the coefficient
 $e^{i (\theta_1+\theta_2)}$ occurs frequently. we can expand and simplify the determinant as follow:
 \begin{small}
 \begin{eqnarray*}
 & &\Phi(\rho)\\
 &=&-(C_1'+e^{i \theta _2}C_2'+e^{i \theta _1}C_3'+e^{i (\theta_1+\theta_2)}C_4')e^{i (\theta_1+\theta_2)}
 [S_1S_2S_3+S_4(e^{i \theta _2}S_1S_2+e^{i \theta _1}S_1S_3+e^{i (\theta_1+\theta_2)}S_2S_3)]\\
 & &+(C_1-e^{i (\theta_1+\theta_2)}C_4)e^{i (\theta_1+\theta_2)}[S_1'S_2S_3+S_1S_2'S_3+S_1S_2S_3'+S_4'(e^{i \theta _1}S_1S_3+e^{i \theta _2}S_1S_2+e^{i
 (\theta_1+\theta_2)}S_2S_3)]\\
 & &+(C_1-e^{i \theta_1}C_3)e^{i(\theta_1+2\theta_2)}[S_1'S_2S_4+S_1S_2'S_4+S_1S_2S_4'+S_3'(e^{i \theta _1}S_2S_4+e^{-i \theta _2}S_1S_2+e^{i
 (\theta_2-\theta_1)}S_1S_4)]\\
 & &-(C_1-e^{i \theta_2}C_2)e^{i(2\theta_1+\theta_2)}[S_1'S_3S_4+S_1S_3'S_4+S_1S_3S_4'+S_2'(e^{-i \theta _1}S_1S_3+e^{i \theta _2}S_3S_4+e^{i
 (\theta_1-\theta_2)}S_1S_4)]\\
 &=&-e^{2i (\theta_1+\theta_2)}[(C_1S_2S_3S_4)'+(S_1C_2S_3S_4)'+(S_1S_2C_3S_4)'+(S_1S_2S_3C_4)']\\
 & &+(C_1S_1'-C_1'S_1)(S_2S_3e^{i(\theta_1+\theta_2)}+S_2S_4e^{i(\theta_1+2\theta_2)}+S_3S_4e^{i(2\theta_1+\theta_2)})\\
 & &+(C_2S_2'-C_2'S_2)(S_1S_3e^{i(\theta_1+2\theta_2)}+S_1S_4e^{i(3\theta_1+\theta_2)}+S_3S_4e^{i(2\theta_1+3\theta_2)})\\
 & &+(C_3S_3'-C_3'S_3)(S_1S_2e^{i(2\theta_1+\theta_2)}+S_1S_4e^{i(\theta_1+3\theta_2)}+S_2S_4e^{i(3\theta_1+2\theta_2)})\\
 & &+(C_4S_4'-C_4'S_4)(S_1S_2e^{i(2\theta_1+3\theta_2)}+S_1S_3e^{i(3\theta_1+2\theta_2)}+S_2S_3e^{i(3\theta_1+3\theta_2)})\\
 &=&e^{2i (\theta_1+\theta_2)}\{2(S_1S_3+S_2S_4)\cos\theta_1+2(S_1S_2+S_3S_4)\cos\theta_2+2S_2S_3\cos(\theta_1+\theta_2)\\
 & &+2S_1S_4\cos(\theta_1-\theta_2)-(C_1S_2S_3S_4)'-(S_1C_2S_3S_4)'-(S_1S_2C_3S_4)'-(S_1S_2S_3C_4)'\}
 \end{eqnarray*}
 \end{small}
 Note that by Lagrange identity, $C_jS'_j-S_jC'_j=1,\ j=1,\ 2,\ 3,\ 4.$ Thus we have proved the following Theorem.
 \newtheorem{th3.2}{Theorem}[section]
 \begin{th3.2}
 \label{th3.2}
 The dispersion relation for the Schr\"{o}dinger operator on the rectangular graph $\Gamma$ as (\ref{eq3.1})
 and (\ref{eq3.4}) are given by with $\theta_j\in[-\pi, \pi],\ j=1,\ 2$
 \begin{eqnarray*}
 & &(C_1S_2S_3S_4)'+(S_1C_2S_3S_4)'+(S_1S_2C_3S_4)'+(S_1S_2S_3C_4)'\\
 &=&2[(S_1S_3+S_2S_4)\cos\theta_1+(S_1S_2+S_3S_4)\cos\theta_2+S_2S_3\cos(\theta_1+\theta_2)+S_1S_4\cos(\theta_1-\theta_2)]
 \end{eqnarray*}
 \end{th3.2}

  If $a_1=a_2$ and $q_j's$ are identical to $q$, then the dispersion relation can be simplified to
 \begin{eqnarray*}
 4S^2(C'S+3CS')&=&2S^2(2\cos\theta_1+2\cos\theta_2+\cos(\theta_1+\theta_2)+\cos(\theta_1-\theta_2))\\
               &=&4S^2(\cos\theta_1+\cos\theta_2+\cos\theta_1\cos\theta_2)\\
               &=&4S^2\, (4\cos^2\frac{\theta_1}2\cos^2\frac{\theta_2}2-1),
 \end{eqnarray*}
  after an elementary computation.
  This implies
  $$\displaystyle4S^2[CS'-\cos^2\frac{\theta_1}2\cos^2\frac{\theta_2}2]=0.
  $$
  Moreover, if $q$ is even, then  $C(a,\rho)=S'(a,\rho)$. Therefore $\rho$ is
 a singular value if and only if $S(a, \rho)=0$ or
 $$S'(a,\rho)=\displaystyle\pm\cos\frac{\theta_1}2\cos\frac{\theta_2}2
 \quad\in [-1,1].$$

 \setcounter{equation}{0}
 \section{Dispersion relation for periodic $n$-cubic graphs}
 \hskip0.25in Consider $Hy(x)=-y''(x)+q(x)y(x)=\lambda y(x)$ where $x\in\mathbb{R}$, $q(x)\in L^2(\mathbb{R})$. For $n$-cubic
 graph with length $a_1,\ a_2,\ \ldots,\ a_n$, we can define the fundamental domain $W$ centered at $\bf{0}$ and surrounded
 by $2n$ hyperplanes.

 We denoted the set of vertices of $W$ as $V_{W}$:
 $$V_{W}:=\{(\pm a_1,0, \ldots, 0), (0, \pm a_2, 0, \ldots, 0), \ldots, (0,\ldots, 0, \pm a_n),
 (\mp a_1, \pm a_2, \ldots, \pm a_{n-1}, \pm a_n)\in\mathbb{R}^n\}$$
 It is clear that the cardinality of $V_{W}=2n+2$ and from $V_{W}$ we can define $n$ vectors to generalize $W$.
 Without loss of generality, these $n$ vectors are started at $v_0=(-a_1, 0, \ldots, 0)$ and denoted as $k_j, j=1, \ldots, n$.
 \begin{eqnarray*}
 \left\{
 \begin{array}{lll}
   k_1 &=& (0, a_2, a_3, \ldots, a_n)\\
   k_2 &=& (a_1, -a_2, 0, \ldots, 0)\\
   \vdots\\
   k_n &=& (a_1, 0, \ldots, 0, -a_n)
 \end{array}
 \right.
 \end{eqnarray*}
 At $v_0$, $\{v_0+k_2,\ \ldots,\ v_0+k_n\}$ forms the hyperplane $\displaystyle \frac{x_1}{a_1}+\ldots+\frac{x_n}{a_n}=-1$, $\{v_0+k_1,\ \ldots,\ v_0+k_{n-1}\}$
 forms the hyperplane $\displaystyle \frac{x_1}{a_1}+\ldots+\frac{x_{n-1}}{a_{n-1}}-(n-2)\frac{x_n}{a_n}=-1$, $\ldots$,
 $\{v_0+k_1,\ v_0+k_3,\ \ldots,\ v_0+k_n\}$ forms the hyperplane $\displaystyle\frac{x_1}{a_1}-(n-2)\frac{x_2}{a_2}+\ldots+\frac{x_n}{a_n}=-1$.
 Thus, the fundamental domain $\Gamma_0$ is defined as (see \eqref{fig3.2} as an example)
 \begin{small}
 \begin{eqnarray*}
 \Gamma_0&:=&\Bigl\{(x_1, x_2, \ldots, x_n):\displaystyle -1\leq\sum^n_{j=1}\frac{x_j}{a_j}\leq (n-2)\\
 & &\mbox{ and }-1\leq\sum^n_{j=1}\frac{x_j}{a_j}-(n-1)\frac{x_k}{a_k}\leq (n-2),\ k=2, \ldots, n\Bigr\}.
 \end{eqnarray*}
 \end{small}
 \begin{figure}[h!]
 \centering\includegraphics[width=15cm,height=10cm]{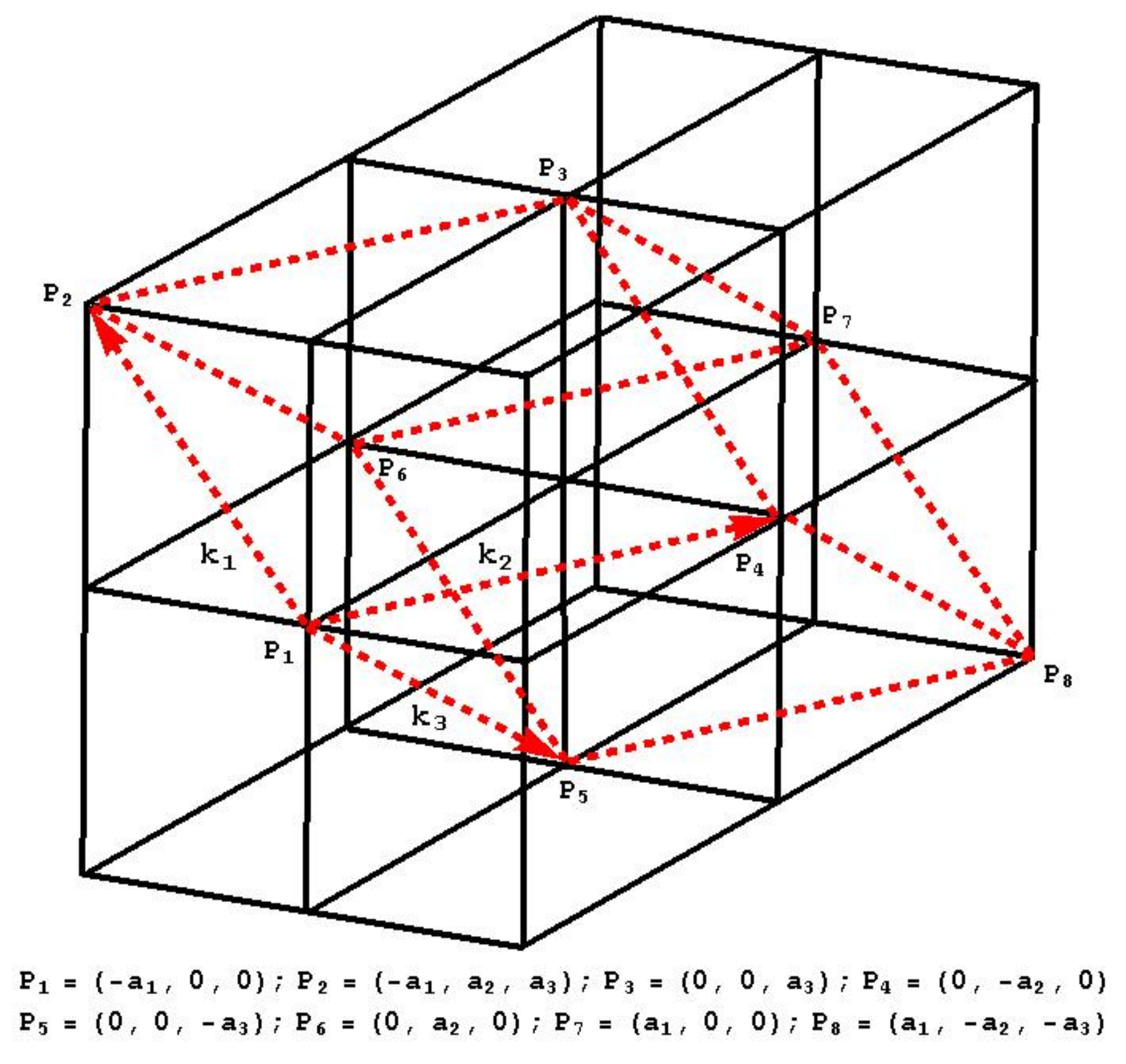}
 \caption{Fundamental domain for 3-cubic graph}
 \label{fig3.3}
 \end{figure}
 $W$ contains $2n$ edges of the $n$-cubic graph and on each edge we have the equation:
 \begin{equation}
 Hy_j''(x)=-y_j''(x)+q_j(x)y_j(x)=\lambda y_j(x),\ j=1, \ldots, 2n. \label{eq3.10}
 \end{equation}

 By the Floquet-Bloch conditions, for any solution of (\ref{eq3.10}), we have
 \begin{equation*}
 \displaystyle y_j(x+\sum^n_{l=1}p_lk_l)=e^{i\sum^n_{l=1}p_l\theta_l}y_j(x),\ \theta_l\in[-\pi, \pi].
 \end{equation*}
 So, we have the following equations:
 \begin{equation}
 \left\{
 \begin{array}{lll}
 y_1(0, \lambda)=y_2(0, \lambda)=\ldots=y_{2n}(0, \lambda)\\
 \displaystyle\sum^{2n}_{j=1}y'_j(0, \lambda)=0\\
 y_1(a_1,\lambda)=e^{i \theta_j}y_j(a_j,\lambda),\ (j=2,\ \ldots,\ n)\\
         \qquad=e^{i(\Theta-\theta_{2n-j+1})}y_{j}(a_{2n-j+1},\lambda),\ (j=n+1,\ \ldots,\ 2n-1)=e^{i\Theta}y_{2n}(a_1, \lambda)\\
 \displaystyle y_1'(a_1,\lambda)+\sum^n_{j=2}e^{i \theta_j}y_j'(a_j,\lambda)+\sum^{2n-1}_{k=n+1}e^{i(\Theta-\theta_{2n-k+1})}y_{k}'(a_{2n-j+1},\lambda)+e^{i
 \Theta}y_{2n}(a_1,\lambda)=0,
 \end{array}
 \right. \label{eq3.11}
 \end{equation}
 where $\Theta=\displaystyle\sum^n_{j=1}\theta_j$.

 Let $\rho=\sqrt{\lambda}$, then any solution $y_j$  satisfies
 $$
 y_j(x, \rho)=A_j(\rho)C_j(x, \rho)+B_j(\rho)S_j(x, \rho),\ j=1,\ldots, 2n.
 $$
 If we write $C_j=C_j(a,\rho)$, $S_j=S_j(a,\rho)$, $A_j=A_j(\rho)$ and $B_j=B_j(\rho)$ for short, then
 $$
 \left\{
 \begin{array}{l}
 A_1=A_2=\ldots=A_{2n}\\
 \displaystyle\sum^{2n}_{j=1}B_j=0\\
 A_1C_1+B_1S_1=e^{i \theta_j}(A_jC_j+B_jS_j),\ (j=2,\ \ldots,\ n)\\
       \qquad=e^{i(\Theta-\theta_{2n-j+1})}(A_jC_j+B_jS_j),\ (j=n+1,\ \ldots,\ 2n-1)=e^{i\Theta}(A_{2n}C_{2n}+B_{2n}S_{2n})\\
 (A_1C'_1+B_1S'_1)+\displaystyle\sum^{n}_{j=2}e^{i \theta_j}(A_jC'_j+B_jS'_j)+\sum^{2n}_{k=n+1}e^{i (\Theta-\theta_{2n-k+1})}(A_kC_k'+B_kS'_k)\\
       \qquad+e^{i \Theta}(A_{2n}C'_{2n}+B_{2n}S'_{2n})=0
 \end{array}
 \right.
 $$
 \vskip0.2in
 Define
 \begin{equation}
 \tau_j=\left\{\begin{array}{llll}
 1&,\ j=1\\
 \displaystyle e^{i\theta_j}&,\ 2\leq j\leq n \\
 \displaystyle e^{i(\Theta-\theta_{2n-j+1})}&,\ n+1\leq j\leq2n-1 \\
 e^{i\Theta}&,\ j=2n.
 \end{array}\right. \label{eq3.35}
 \end{equation}

 Though $\theta_1$ is arbitrary between $[-\pi, \pi]$, while $\theta_1=0$, we may reduce $\tau_j$ as:
 $$
 \tau_j=\left\{\begin{array}{llll}
 \displaystyle e^{i\theta_j}&,\ 1\leq j\leq n \\
 \displaystyle e^{i(\Theta-\theta_{2n-j+1})}&,\ n+1\leq j\leq2n. \\
 \end{array}\right.
 $$

 It is clear that the characteristic function $\Phi(\rho)$ satisfies
 \begin{eqnarray*}
 \Phi(\rho)=
 \left|
 \begin{array}{ccccccccc}
 0 & 1 & 1 & \cdots & \cdots & \cdots & 1 \\
 C_1-\tau_2C_2 & S_1 & -\tau_2S_2 & 0 & \cdots & \cdots & 0 \\
 C_1-\tau_3C_3 & S_1 & 0 & -\tau_3S_3 & \ddots & \cdots & 0 \\
 \vdots & \vdots & \vdots & \ddots & \ddots & \ddots & \vdots\\
 C_1-\tau_{2n-1}C_{2n-1} & S_1 & 0 & \cdots & \cdots &-\tau_{2n-1}S_{2n-1} & 0 \\
 C_1-\tau_{2n}C_{2n} & S_1 & 0 & \cdots & \cdots & 0 & -\tau_{2n}S_{2n}\\
 \displaystyle\sum^{2n}_{k=1} \tau_k C'_k & S_1' & \tau_2S'_2 & \tau_3S_3' & \cdots & \tau_{2n-1}S_{2n-1}' & \tau_{2n}S_{2n}'\\
 \end{array}
 \right|
 \end{eqnarray*}

 \newtheorem{th3.4}{Theorem}[section]
 \begin{th3.4}
 \label{th3.4}
 For $n$-cubic quantum graph, let $\d \Th=\sum_{1}^n \th_j$, and $\tau_j$ be given in \eqref{eq3.35}.
 \begin{enumerate}
 \item[(a)] The characteristic function
 $$
 \Phi(\rho)=-(\sum^{2n}_{k=1}\tau_kC_k')\sum^{2n}_{l=1}(\prod_{i\neq l}\tau_iS_i)
 +\sum^{2n}_{l=2}\left((C_1-\tau_lC_l)\sum_{k\neq l}\frac{\prod_{i\neq l,k}(\tau_iS_i)}{\tau_kS_k}(\tau_kS_k'-\tau_lS_l')\right).$$
 \item[(b)] If all the edgelengths $a_i$'s are equal to $a$, and $q_i$'s are identical and even on $(0, a)$,  then $C_j=C$ and $S_j=S$, for all $j$,
 $$
 \displaystyle\Phi(\rho)=4e^{in\Theta}S^{2n-2}\left\{\left[\cos(\frac{\Theta}2)
                        +\sum^n_{l=2}\cos(\frac{\Theta-2\theta_l}2)\right]^2-n^2 \, (S')^2\right\}.
 $$
 Thus $\rho$ is a singular value if and only if
 \begin{equation}
 \displaystyle S(a, \rho)=0\ \mbox{or }S'(a,\rho)=\pm\frac 1 n
 \left(\cos(\frac{\Theta}2)+\sum^n_{l=2}\cos(\frac{\Theta-2\theta_l}2)\right).\label{eq3.5}
 \end{equation}
 \end{enumerate}
 \end{th3.4}
   The determinant $\Phi(\rho)$ above involves a $n\times n$ matrix of a special form.  We decompose it into two parts and compute their values
   accordingly.
 \newtheorem{lem3.5}[th3.4]{Lemma}
 \begin{lem3.5}
 \label{lem3.5}
 For $m\in\mathbb{N}$, let $a_1,\ldots,a_m\in\mathbb{C}$. Then the determinants $\{u_m\},\ \{v_{m,l}\}$ defined below are
 given as:
 \begin{footnotesize}
 \begin{enumerate}
 \item[(a)] $u_m:=\left|
             \begin{array}{ccccccccccc}
             1 & \cdots & \cdots & \cdots & \cdots & 1\\
             a_1 & -a_2 & 0 & \cdots & \cdots & 0\\
             \vdots & 0 & \ddots & \ddots & \cdots & \vdots \\
             \vdots & \vdots & \ddots & \ddots & \ddots & \vdots\\
             \vdots & \vdots & \cdots & \ddots & \ddots & 0 \\
             a_1 & 0 & \cdots & \cdots & 0 & -a_m
             \end{array}
             \right|
             =(-1)^{m-1}\displaystyle\sum^m_{j=1}(\prod_{i\neq j}a_i)$
 \item[(b)] $
             v_{m,l}:=\left|
             \begin{array}{ccccccccccc}
             1 & 1 & \cdots & \cdots & 1 & \cdots & \cdots &1\\
             a_1 & -a_2 & 0 & \cdots & \cdots & \cdots & \cdots & 0\\
             \vdots & 0 & \ddots & \ddots & \cdots & \cdots & \cdots & 0\\
             \vdots & \vdots & \ddots & -a_{l-1} & 0 & \cdots & \cdots & 0\\
             \vdots & \vdots & \cdots & \ddots & 0 & -a_{l+1} & \ddots & 0\\
             \vdots & \vdots & \cdots & \cdots & \ddots & \ddots & \ddots & \vdots\\
             a_1 & 0 & \cdots & \cdots & \cdots & \cdots & 0 &-a_m\\
             b_1 & b_2 & \cdots & b_{l-1} & b_l & b_{l-1} & \cdots & b_m \\
             \end{array}
             \right|
             =(-1)^{l+1}\displaystyle\sum^n_{\overset{k=1}{k\neq l}}\frac{\prod_{i\neq l} a_i}{a_k}(b_k-b_l),
             $
              $l=2,\cdots,\ m.$
 \end{enumerate}
 \end{footnotesize}
 \end{lem3.5}
 \begin{proof}
 \begin{enumerate}
 \item[(a)] The determinant is easily computed by a cofactor expansion along the first row.
 \item[(b)] Let us perform cofactor expansion along the $l^{th}$ column. Then
            \begin{small}
            \begin{eqnarray*}
            v_{m,l}
            &=&(-1)^{l+1}\left|
            \begin{array}{cccccccc}
            a_1 & -a_2 & 0 & \cdots & \cdots & \cdots & 0\\
            \vdots & 0 & \ddots & \ddots & \cdots & \cdots & \vdots\\
            \vdots & \vdots & \ddots & -a_{l-1} & \ddots & \cdots & \vdots\\
            \vdots & \vdots &\cdots & \ddots & -a_{l+1} & \ddots & \vdots\\
            \vdots & \vdots & \cdots & \cdots & \ddots & \ddots &\vdots\\
            a_1 & 0 & \cdots & \cdots & \cdots & 0 & -a_m\\
            b_1 & \cdots & \cdots & b_{l-1} & b_{l+1} & \cdots & b_m
            \end{array}
            \right|\\
            &+&(-1)^{l+m}b_l\left|
            \begin{array}{cccccccc}
            1 & 1 & \cdots & 1 & 1 & \cdots & 1\\
            a_1 & -a_2 & 0 & \cdots & \cdots & \cdots & 0\\
            \vdots & 0 & \ddots & \ddots & \cdots & \cdots & \vdots\\
            \vdots & \vdots &\ddots & -a_{l-1} & \ddots & \cdots & \vdots\\
            \vdots & \vdots & \ddots & \ddots & -a_{l+1} & \ddots & \vdots\\
            \vdots & \vdots & \cdots & \cdots & \ddots & \ddots & \vdots\\
            a_1 & 0 & \cdots & \cdots & \cdots & 0 & -a_m
            \end{array}
            \right|.
            \end{eqnarray*}
            \end{small}
            By (a), the second term equals $\displaystyle(-1)^{l+2m-2}b_l\sum^n_{\overset{k=1}{k\neq l}}(\frac{\prod_{i\neq l}a_i}{a_k})$.
            Then we perform cofactor expansion along the last row to see that the first term of $v_{m,l}$
            is exactly $\displaystyle(-1)^{l+2m+1}\sum^n_{\overset{k=1}{k\neq l}}b_k(\frac{\prod_{i\neq l}a_i}{a_k})$.
            The proof is complete.
            \end{enumerate}
 \end{proof}

 \begin{proof}[Proof of Theorem \ref{th3.4}]
 (a) \ Using cofactor expansion along the first column, we have
 \begin{eqnarray*}
 \Phi(\rho)&=&-(C_1-\tau_2C_2)
 \left|
 \begin{array}{ccccccccc}
 1 & 1 & \cdots & \cdots & \cdots & 1 \\
 S_1 & 0 & -\tau_3S_3 & 0 & \cdots & 0 \\
 \vdots & \vdots & \ddots & \ddots & \ddots & \vdots\\
 S_1 & 0 & \cdots & 0 &-\tau_{2n-1}S_{2n-1} & 0 \\
 S_1 & 0 & \cdots & \cdots & 0 & -\tau_{2n}S_{2n}\\
 S_1' & \tau_2S'_2 & \tau_3S_3' & \cdots & \tau_{2n-1}S_{2n-1}' & \tau_{2n}S_{2n}'\\
 \end{array}
 \right|+\cdots\\
 & &+\displaystyle(\sum^{2n}_{k=1} \tau_k C'_k)
 \left|
 \begin{array}{ccccccccc}
 1 & 1 & \cdots & \cdots & \cdots & 1 \\
 S_1 & -\tau_2S_2 & 0 & \cdots & \cdots & 0 \\
 S_1 & 0 & -\tau_3S_3 & 0 & \cdots & 0 \\
 \vdots & \vdots & \ddots & \ddots & \ddots & \vdots\\
 S_1 & 0 & \cdots & 0 &-\tau_{2n-1}S_{2n-1} & 0 \\
 S_1 & 0 & \cdots & \cdots & 0 & -\tau_{2n}S_{2n}\\
 \end{array}
 \right|.
 \end{eqnarray*}
 Set $a_i=\tau_iS_i$, and $b_i=\tau_iS_i'$ , then by Lemma \ref{lem3.5}, expanding each term and reducing the formula,
 we obtain the result.
 \\[0.1in]
 (b) \ From (a), when $C_j=C$ and $S_j=S$,
 \begin{eqnarray*}
 \displaystyle\Phi(\rho)&=&-C'S^{2n-1}(\sum^{2n}_{k=1}\tau_k)\sum^{2n}_{l=1}(\prod_{i\neq l,k}\tau_i)
 +CS^{2n-2}S'\sum^{2n}_{l=2}(1-\tau_l)\sum_{k\neq l}\frac{\prod_{i\neq l}(\tau_i)}{\tau_k}(\tau_k-\tau_l)\\
 &=&(\prod^{2n}_{i=1}\tau_i)\left(-C'S^{2n-1}(\sum^{2n}_{k=1}\tau_k)(\sum^{2n}_{l=1}\frac1{\tau_l})
 +CS^{2n-2}S'\sum^{2n}_{l=2}\left[(1-\tau_l)\sum_{k\neq l}(\frac1{\tau_l}-\frac1{\tau_k})\right]\right)\\
 &=&e^{in\Theta}\, \left\{-C'S^{2n-1}(\sum^{2n}_{k=1}\tau_k)(\sum^{2n}_{l=1}\frac1{\tau_l})
 +CS^{2n-2}S'\sum^{2n}_{l=2}\left[(1-\tau_l)(\frac{2n}{\tau_l}-\sum^{2n}_{k=1}\frac{1}{\tau_k})\right]\right\}\\
 &:=&e^{in\Theta}\, (-I_1+I_2)
 \end{eqnarray*}
  Now by the Euler identity,
 \begin{eqnarray*}
 \sum^{2n}_{k=1}\tau_k\cdot\sum^{2n}_{l=1}\frac1{\tau_l}
 &=&\rme^{i \Theta}\, \left(1+\sum_{j=2}^n (\rme^{-i \theta_j}+ \rme^{-i (\Theta-\theta_j)})+\rme^{-i \Theta}\right)^2\\
 &=&\displaystyle \left(\rme^{i \frac{\Theta}{2}}+\sum_{j=2}^n(\rme^{i (\frac{\Theta-2\theta_j}{2})}+\rme^{-i (\frac{\Theta-2\theta_j}{2})})+\rme^{-i \frac{\Theta}{2}}\right)^2\\
 &=&4\displaystyle\left(\cos(\frac{\Theta}2)+\sum^n_{l=2}\cos(\frac{\Theta-2\theta_l}2)\right)^2\\
 &:=&4T.
 \end{eqnarray*}
 Hence $\d
 I_1= 4T C'S^{2n-1}$.
 On the other hand,
 \begin{eqnarray*}
 I_2&=&CS^{2n-2}S'\displaystyle\sum^{2n}_{l=2}\left[(1-\tau_l)(\frac{2n}{\tau_l}-\sum^{2n}_{k=1}\frac1{\tau_k})\right]\\
 &=&CS^{2n-2}S'\left[\displaystyle 2n\sum^{2n}_{l=2}\frac{1}{\tau_l}-2n(2n-1)-(2n-1)\sum^{2n}_{k=1}\frac1{\tau_k}+
 \sum^{2n}_{l=1}\tau_l\cdot\sum^{2n}_{k=1}\frac1{\tau_k}-\sum^{2n}_{k=1}\frac1{\tau_k}\right]\\
 &=&4(T-n^2)CS^{2n-2}S'.
 \end{eqnarray*}

 So by above and the Lagrange identity, we have
 \begin{eqnarray*}
 \Phi(\rho)&=&4\{T-n^2CS'\}\rme^{in\Theta}S^{2n-2}\\
           &=&4\left\{\left(\cos(\frac{\Theta}2)+\sum^n_{l=2}\cos(\frac{\Theta-2\theta_l}2)\right)^2-n^2CS'\right\}\rme^{in\Theta}S^{2n-2}.
 \end{eqnarray*}

 If $q$ is even, then $C=S'$, $\rho$ is a singular value if and only if
 \begin{eqnarray*}
 \Phi(\rho)&=&4\left\{\left(\cos(\frac{\Theta}2)+\sum^n_{l=2}\cos(\frac{\Theta-2\theta_l}2)\right)^2-n^2S'^2\right\}e^{in\Theta_n}S^{2n-2}=0.
 \end{eqnarray*}
 That is , if and only if
 $$\displaystyle S(a, \rho)=0\ \mbox{or }S'(a,\rho)=\pm\frac 1 n
 \left(\cos(\frac{\Theta}{2})+\sum^n_{l=2}\cos(\frac{\Theta-2\theta_l}2)\right).$$
 \end{proof}

 Thus the spectrum $\sigma(H)$ is given by
 $$
 \sigma(H)=\left\{\rho^2:\displaystyle S(a, \rho)=0\ \mbox{or }S'(a,\rho)=\pm\frac 1 n
 \left(\cos(\frac{\Theta}{2})+\sum^n_{l=2}\cos(\frac{\Theta-2\theta_l}2)\right)\right\}.
 $$
 It can be easily seen that the case $S(a, \rho)=0$ gives eigenvalues of infinite multiplicities.
 \newtheorem{th3.6}[th3.4]{Theorem}
 \begin{th3.6}
 \label{th3.6}
 For periodic $n$-cubic graph when $a_i$'s are the same, and $q_i$'s are identical and even,
 \begin{equation}
  \sigma_{p}(H)=\{\lambda=\rho^2:\ S(a,\rho)=0\} \label{eq3.12}
 \end{equation}
  and each $\lambda\in \sigma_p(H)$ has infinite multiplicity.
 \end{th3.6}

 \begin{proof}

 \begin{figure}[h!]
 \centering\includegraphics[width=7.1cm,height=5cm]{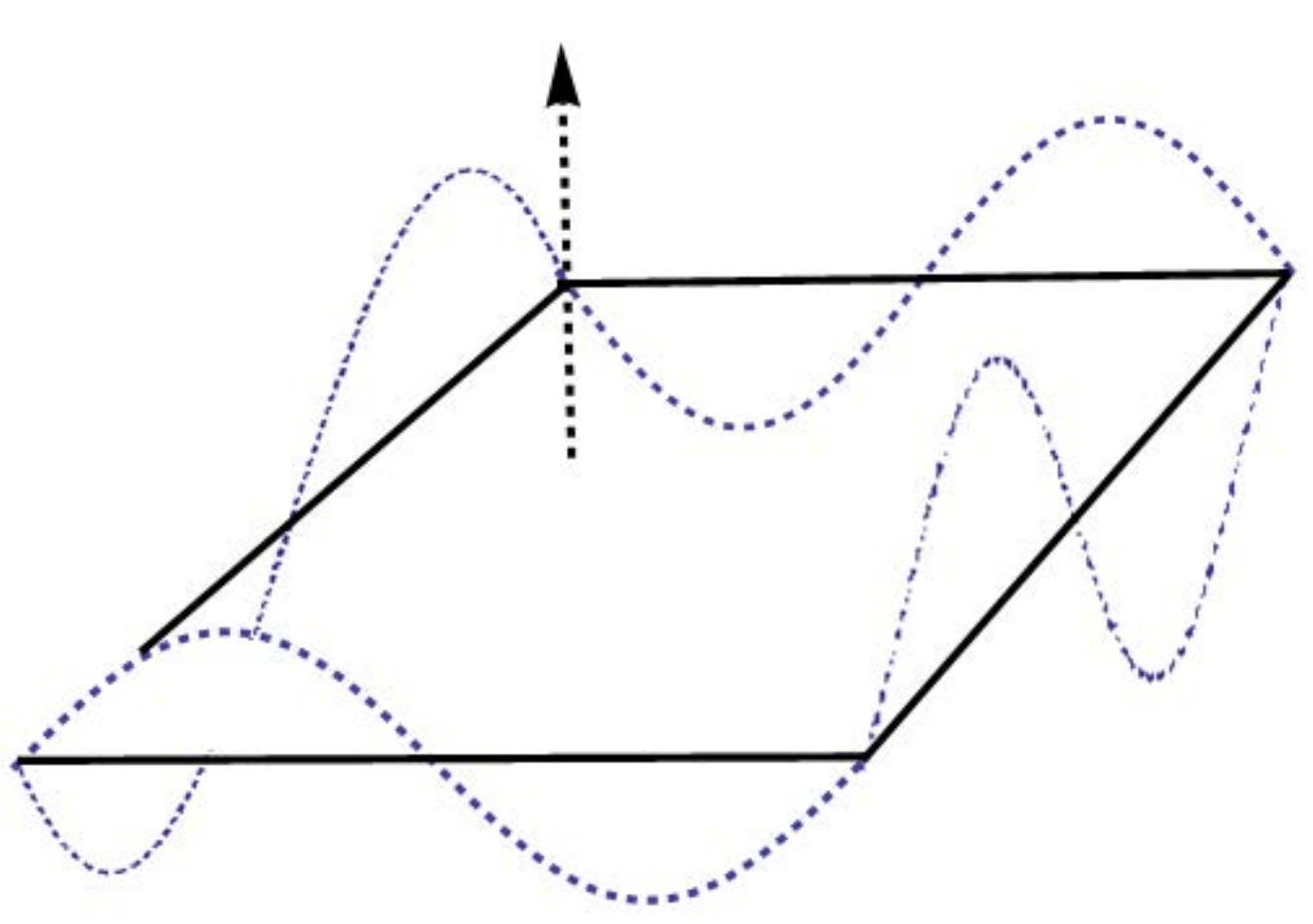}
 \caption{Characteristic function for $\mathcal{S}$}
 \label{fig2.5}
 \end{figure}

 Let $\rho$ satisfy $S(a,\rho)=0$. Then there is a Dirichlet eigenfunction $f$ of the Schr\"{o}dinger operator on the interval $(0,a)$.
 Reproduce this function $f$ on the edges of a square in the graph (see Figure \ref{fig2.5}). Then extend the result to $\widetilde{f}$ on the
 whole graph by letting its value be $0$ outside the square. In this way, $\widetilde{f}\in L^2(\Gamma)$ where $\Gamma$ is the periodic graph of squares.
 It also satisfies the vertex conditions and hence
 $$
 H\widetilde{f}=\rho^2\widetilde{f}.
 $$
 That is $\widetilde{f}$ is an eigenfunction. Moreover $\rho^2$ is an eigenvalue of infinite multiplicity, for there is an eigenfunction associated with every
 square in the periodic graph.
 \end{proof}

 Furthermore, the absolutely continuous spectrum $\sigma_{ac}(H)$ is defined in \eqref{eq3.5}. It is easy to see that for any $\la=\rho^2\in \sigma_{ac}(H)$,
 $S'(a,\rho)\subset [-1,1]$. On the other hand, when $\th_j=0$ for all $j$, \eqref{eq3.5} implies that $S'(a,\rho)=\pm 1$. As $\phi=S'(a,\cdot)$ is entire in
 $\rho$, we have proved the following theorem
 \newtheorem{th3.7}[th3.4]{Theorem}
 \begin{th3.7}
 \label{th3.7}
 For the $n$-cubic periodic quantum graph,
 $$
 \sigma_{ac}(H)=\{ \rho^2\in \bfR:\ | S'(a,\rho)\in [-1,1]\} = \phi^{-1}([-1,1]).
 $$
 \end{th3.7}
 \section{Concluding remarks}
  As a summary, we have derived the dispersion relations of the $n$-cubic periodic quantum graphs, for $n=2$ and beyond.  We do not require the edgelengths to be the same.
  Neither are the potentials required to be identical and even.  The spectra are all
  of a band and gap structure. When the edgelengths $a_i$'s are the same, and the potentials $q_i$'s are identical and even, the dispersion relations become simple.
  The point spectrum $\sig_p(H)$ coincide with the Dirichlet-Dirichlet eigenvalues on the interval, and corresponding eigenfunctions has infinite multiplicity.
  Furthermore, the absolutely continuous spectra satisfy $\sigma_{ac}(H) = \phi^{-1}([-1,1])$, where $\phi(\rho)=S'(a,\rho)$.
  As it is well known that as $\rho\to\infty$,
  $$
  S'(a,\rho)-\cos(\rho a)=O(\frac{1}{\rho}).
  $$
  we know that the spectral gaps for the spectrum converges to $0$.  Our analysis is rigorous and general. We expect that our method should be able to
  find the dispersion relations for the energies of many three dimensional crystal lattices.

  \section*{Acknowledgements}
  We thank  Eduardo Jatulan and Vyacheslav Pivovarchik for stimulating discussions.
The authors are partially supported by Ministry of Science and Technology, Taiwan, under contract number MOST105-2115-M-110-004.

\end{document}